\documentclass[11pt]{article}

\usepackage{graphicx}

\def\be{\begin{equation}}
\def\ee{\end{equation}}
\newcommand{\kk}[2]{\frac{#1}{#2}}

\newcommand{\ff}[1]{{\bf  #1}}

\def\lam{\lambda}
\def\s{\qquad}

\def\={\approx}
\def\8{{\infty}}
\def\x{\ff{x}}
\def\vcode#1#2#3#4{\begin{figure}
\begin{center}
\begin{minipage}[c]{#1\textwidth}
{{\small #2 \hrule \vspace{5pt}   %
{\it #3}  \vspace{5pt} \hrule }}
\end{minipage}
\caption{#4}
\end{center}   \end{figure}} %%
\begin{document}

\title{Flower Pollination Algorithm for Global Optimization}

\author{Xin-She Yang \\ \\
Department of Engineering, University of Cambridge, \\
Trumpington Street, Cambridge CB2 1PZ, UK. }

\date{}

\maketitle

\begin{abstract}
Flower pollination is an intriguing process in the natural world. Its evolutionary
characteristics can be used to design new optimization algorithms. In this paper,
we propose a new algorithm, namely, flower pollination algorithm,
inspired by the pollination process of flowers. We first use ten test functions
to validate the new algorithm, and compare its performance with genetic algorithms
and particle swarm optimization. Our simulation results show the flower algorithm
is more efficient than both GA and PSO. We also use the flower algorithm to solve
a nonlinear design benchmark, which shows the convergence rate is almost exponential.
\end{abstract}

{\bf Citation Details:} Xin-She Yang, Flower pollination algorithm for global optimization,
in: {\it Unconventional Computation and Natural Computation 2012}, Lecture Notes in
Computer Science, Vol. 7445, pp. 240-249 (2012).

%% Begin of Main Text %%

\section{Introduction}

Nature has been solving challenging problems over millions and billions of years,
and many biological systems have evolved with intriguing and surprising efficiency
in maximizing their evolutionary objectives such as reproduction. Based on the successfully
characteristics of biological systems, many nature-inspired algorithms have been developed
over the last few decades \cite{Yang,YangBook}. For example, genetic algorithms
were based on the Darwinian evolution of biological systems \cite{Holland} and
particle swarm optimization was based on the swarm behaviour of birds and fish \cite{Ken,Ken2},
which bat algorithm was based on the echolocation behaviour of microbats \cite{YangBat} and
firefly algorithm was based on the flashing light patterns of tropic fireflies \cite{Yang,YangFA}.
All these algorithms have been applied to a wide range of applications.

In many design applications in engineering and industry, we have to try to find the
optimal solution to a given problem under highly complex constraints. Such constrained
optimization problems are often highly nonlinear, to find the optimal solutions is often
a very challenging task if it is not impossible. Most conventional
optimization do not work well for problems with nonlinearity and multimodality.
Current trend is to use nature-inspired metaheuristic algorithms to tackle such
difficult problems, and it has been shown that metaheuristics are surprisingly very efficient.
For this reason, the literature of metaheuristics has expanded tremendously in the last
two decades \cite{Yang,YangBook}. Up to now, researchers have only use a very limited
characateristics inspired by nature, and there is room for more algorithm development.

In this paper, we will propose a new algorithm based on the flower pollination process
of flowering plants. We will first briefly review the main characteristics of
flower pollination, and thus idealize these characteristics into four rules. We will then
use them to develop a flower pollination algorithm (FPA), or the flower algorithm.
Then, we validate it using a set of well-known test functions and design benchmark.
We analyze the simulations and make comparison of its performance with genetic algorithm
and particle swarm optimization. Finally, we discuss further topics for extending this algorithm.

From the biological evolution point of view, the objective of the flower pollination
is the survival of the fittest and the optimal reproduction of plants in terms of
numbers as well as most fittest. This is in fact an optimization process of plant species.
All the above factors and processes of flower pollination interact so as to achieve
optimal reproduction of the flowering plants. Therefore, this can inspire to design
new optimization algorithm. The basic idea of flower pollination in the context
of bees and clustering was investigated before \cite{Kaz}, but in this paper, we will
design a completely new optimization solely based on the flower pollination characteristics.

\section{Characteristics of Flower Pollination}

It is estimated that there are over a quarter of a million types of flowering plants in Nature
and that about 80\% of all plant species are flowering species. It still remains
partly a mystery how flowering plants came to dominate the landscape from
Cretaceous period  \cite{BBCWalker,Flower}.
Flowering plant has been evolving for more than 125 million years
and flowers have become so influential in evolution, we cannot image
how the plant world would be without flowers. The main purpose of a
a flower is ultimately reproduction via pollination. Flower pollination is
typically associated with the transfer of pollen, and such transfer
is often linked with pollinators such as insects, birds, bats and other animals.
In fact, some flowers and insects have co-evolved into a very specialized
flower-pollinator partnership. For example, some flowers can only attract and
can only depend on a specific species of insects for successful pollination \cite{Glover}.

Pollination can take two major forms: abiotic and biotic. About 90\%
of flowering plants belong to biotic pollination, that is, pollen is
transferred by a pollinator such as insects and animals. About 10\%
of pollination takes abiotic form which does not require any pollinators.
Wind and diffusion in water help pollination of such flowering plants
and grass is a good example \cite{Pollen,Glover}.
Pollinators, or sometimes called pollen vectors, can be very diverse.
It is estimate there are at least 200,000 variety of pollinators such as
insects, bats and birds.

Honeybees are a good example of pollinator, and they
can also developed the so-called flower constancy \cite{Chitt}. That is,
these pollinators tend to visit exclusive certain flower species while bypassing
other flower species. Such flower constancy may have evolutionary advantages
because this will maximize the transfer of flower pollen to the same or conspecific plants,
and thus maximizing the reproduction of the same flower species.
Such flower constancy may be advantageous for pollinators as well, because they
can be sure that nectar supply is available with their limited memory and minimum
cost of learning or exploring. Rather than focusing on some unpredictable but potentially
more rewarding new flower species, flower constancy may require minimum investment
cost and more likely guaranteed intake of nectar \cite{Waser}.

Pollination can be achieved by self-pollination or cross-pollination.
Cross-pollination, or allogamy, means pollination can occur from pollen
of a flower of a different plant, while self-pollination is the fertilization
of one flower, such as peach flowers, from pollen of the same flower or different flowers of the same plant,
which often occurs when there is no reliable pollinator available.

Biotic, cross-pollination may occur at long distance, and the pollinators such as
bees, bats, birds and flies can fly a long distance, thus they can considered as the
global pollination. In addition, bees and birds may behave as L\'evy flight behaviour \cite{Pav},
with jump or fly distance steps obey a L\'evy distribution. Furthermore,
flower constancy can be used an increment step using the similarity or difference
of two flowers.

\section{Flower Pollination Algorithm}

Now we can idealize the
above characteristics of pollination process, flower constancy and pollinator behaviour as the
following rules:

\begin{enumerate}

\item Biotic and cross-pollination is considered as global pollination process with
pollen-carrying pollinators performing L\'evy flights.

\item Abiotic and self-pollination are considered as local pollination.

\item Flower constancy can be considered as the reproduction probability is
proportional to the similarity of two flowers involved.

\item Local pollination and global pollination is controlled by a
switch probability $p \in [0,1]$. Due to the physical proximity
and other factors such as wind, local pollination can have a significant
fraction $p$ in the overall pollination activities.

\end{enumerate}

Obviously, in reality, each plant can have multiple flowers, and each flower
patch often release millions and even billions of pollen gametes. However,
for simplicity, we also assume that each plant only has one flower, and each
flower only produce one pollen gamete. Thus, there is no need to
distinguish a pollen gamete, a flower, a plant or solution to a problem.
This simplicity means a solution $\x_i$ is equivalent to a flower and/or a pollen gamete.
In future studies, we can easily extend to multiple pollen gametes for each
flower and multiple flowers for multiobjective optimization problems.

From the above discussions and the idealized characteristics, we can design a flower-based on
algorithm, namely, flower pollination algorithm (FPA). There are two key steps in this algorithm,
they are global pollination and local pollination.

In the global pollination step, flower pollens are carried by pollinators such as insects,
and pollens can travel over a long distance because insects can often fly and move in a much
longer range. This ensures the pollination and reproduction of the most fittest, and thus
we represent the most fittest as $\ff{g}_*$. The first rule plus flower constancy can be  represented mathematically as
\be \x_i^{t+1}=\x_i^t + L (\x_i^t -\ff{g}_*), \ee
where $\x_i^t$ is the pollen $i$ or solution vector $\x_i$ at iteration $t$, and $\ff{g}_*$ is
the current best solution found among all solutions at the current generation/iteration.
The parameter $L$ is the strength of the pollination, which essentially is a step size.
Since insects may move over a long distance with various distance steps, we can use a L\'evy flight
to mimic this characteristic efficiently \cite{Pav,Reynolds}. That is, we draw $L>0$ from a Levy distribution
 \be L \sim \frac{\lam \Gamma(\lam) \sin (\pi \lam/2)}{\pi} \frac{1}{s^{1+\lam}}, \quad (s \gg s_0>0). \ee
Here $\Gamma(\lam$) is the standard gamma function, and this distribution is valid for
large steps $s>0$. In all our simulations below, we have used $\lam=1.5$.

The local pollination (Rule 2) and flower constancy can be represented as
\be \x_i^{t+1}=\x_i^t + \epsilon (\x_j^t -\x_k^t), \ee
where $\x_j^t$ and $\x_k^t$ are pollens from the different flowers of the same plant species.
This essentially mimic the flower constancy in a limited neighborhood. Mathematically,
if $\x_j^t$ and $\x_k^t$ comes from the same species
or selected from the same population, this become a local random walk if we draw
$\epsilon$ from a uniform distribution in [0,1].

\vcode{0.97}{{\sf Flower Pollination Algorithm (or simply Flower Algorithm) }} {
\indent Objective $\min$ or $\max f(\x)$, $\x=(x_1,x_2,..., x_d)$ \\
\indent Initialize a population of $n$ flowers/pollen gametes with random solutions  \\
\indent Find the best solution $\ff{g}_*$ in the initial population \\
\indent Define a switch probability $p \in [0,1]$ \\
\indent  {\bf while} ($t<$MaxGeneration) \\
\indent \qquad {\bf for} $i=1:n$ (all $n$ flowers in the population) \\
\indent \quad \qquad {\bf if} rand $<p$, \\
\indent \qquad \qquad Draw a ($d$-dimensional) step vector $L$ which obeys a L\'evy distribution \\
\indent \qquad \qquad Global pollination via $\x_i^{t+1}=\x_i^t + L (\ff{g}_*-\x_i^t)$ \\
\indent \qquad \quad {\bf else} \\
\indent \qquad \qquad Draw $\epsilon$ from a uniform distribution in [0,1] \\
\indent \qquad \qquad Randomly choose $j$ and $k$ among all the solutions \\
\indent \qquad  \qquad Do local pollination via $\x_i^{t+1}=\x_i^t+\epsilon (\x_j^t-\x_k^t)$ \\
\indent \qquad \quad {\bf end if } \\
\indent \qquad \quad Evaluate new solutions \\
\indent \qquad \quad If new solutions are better, update them in the population \\
\indent \qquad {\bf end for} \\
\indent \qquad \quad Find the current best solution $\ff{g}_*$ \\
\indent  {\bf end while} }
{Pseudo code of the proposed Flower Pollination Algorithm (FPA). \label{fpa-code} }

Most flower pollination activities can occur at both local and global scale.
In practice, adjacent flower patches or flowers in the not-so-far-away neighborhood
are more likely to be pollinated by local flower pollens that those far away.
For this, we use a switch probability (Rule 4) or proximity probability $p$ to switch
between common global pollination to intensive local pollination. To start with,
we can use $p=0.5$ as an initially value and then do a parametric study to find the most appropriate parameter
range. From our simulations, we found that $p=0.8$ works better for most applications.

The above two key steps plus the switch condition can be summarized in the pseudo code shown in Fig. \ref{fpa-code}.

\section{Numerical Results}

Any new optimization should be extensively validated and comparison with other algorithms.
There are many test functions, at least over a hundred well-know test functions
However, there is no agreed set of test functions for validating new algorithms, though
there some review and literature \cite{Ackley,Hedar,YangFA}. In this paper,
we will choose a diverse subset of such test functions to validate our proposed Flower Pollination Algorithm (FPA).

In addition, we will also compare the performance of our algorithm with that of genetic algorithms \cite{Gold} and
particle swarm optimization \cite{Ken,Ken2}. Furthermore, we will also apply FPA to solve a well-known
pressure vessel design benchmark \cite{Cag,Gandomi}.

\subsection{Test Functions}

The Ackley function can be written as
\be f(\x)=-20 \exp\Big[-\kk{1}{5} \sqrt{\kk{1}{d} \sum_{i=1}^d x_i^2}\Big] - \exp\Big[\kk{1}{d} \sum_{i=1}^d \cos (2 \pi x_i)\Big]
+20 +e, \ee
which has a global minimum $f_*=0$ at $(0,0,...,0)$.

The simplest of De Jong's functions is the so-called sphere function
\be f(\x) =\sum_{i=1}^n x_i^2, \s
-5.12 \le x_i \le 5.12, \ee
whose global minimum is obviously $f_*=0$ at $(0,0,...,0)$. This function is unimodal and convex.

Easom's function
\be f(\x)=-\cos(x) \cos(y) \exp \Big[-(x-\pi)^2 + (y-\pi)^2  \Big], \ee
whose global minimum is $f_*=-1$ at $\x_*=(\pi, \pi)$ within $-100 \le x,y \le 100$.
It has many local minima.

Griewangk's function \be f(\x) = \kk{1}{4000} \sum_{i=1}^n x_i^2 - \prod_{i=1}^n \cos (\kk{x_i}{\sqrt{i}}) +1,
\s -600 \le x_i \le 600, \ee
whose global minimum is $f_*=0$ at $\x_*=(0,0,...,0)$. This function is highly multimodal. \\

Michaelwicz's function \be f(\x) = - \sum_{i=1}^n \sin (x_i) \cdot \Big[ \sin ( \kk{i x_i^2}{\pi}) \Big]^{2m}, \ee
where $m=10$, and $0 \le x_i \le \pi$ for $i=1,2,...,n$.
In 2D case, we have
\be f(x,y) = - \sin (x) \sin^{20} (\kk{x^2}{\pi}) - \sin (y) \sin^{20} (\kk{2 y^2}{\pi}), \ee
where $(x,y) \in [0,5] \times [0,5]$. This function has a global minimum
$f_* \= -1.8013$ at $\x_*=(x_*,y_*)=(2.20319, 1.57049)$. \\

Rastrigin's function \be f(\x) = 10 n + \sum_{i=1}^n \Big[ x_i^2 - 10 \cos (2 \pi x_i) \Big],
\s -5.12 \le x_i \le 5.12, \ee
whose global minimum is $f_*=0$ at $(0,0,...,0)$. This function is highly multimodal. \\

Rosenbrock's function \be f(\x) = \sum_{i=1}^{n-1} \Big[ (x_i-1)^2 + 100 (x_{i+1}-x_i^2)^2 \Big], \ee
whose global minimum $f_*=0$ occurs at $\x_*=(1,1,...,1)$ in the domain
$-5 \le x_i \le 5$ where $i=1,2,...,n$. In the 2D case, it is often written as
\be f(x,y)=(x-1)^2 + 100 (y-x^2)^2, \ee
which is often referred to as the banana function.  \\

Schwefel's function \be f(\x) = - \sum_{i=1}^n x_i \sin \Big(\sqrt{|x_i|} \Big),
\s -500 \le x_i \le 500, \ee
whose global minimum $f_*\=-418.9829 n$ occurs at $x_i=420.9687$ where $i=1,2,...,n$. \\

Yang's forest-like function \cite{YangBook}
\be f(\x)=\Big( \sum_{i=1}^d |x_i| \Big) \exp\Big[- \sum_{i=1}^d \sin (x_i^2) \Big],
\;\;\; -2 \pi \le x_i \le 2 \pi, \ee
has a global minimum $f_*=0$ at $(0,0,...,0)$.

Shubert's function \be f(\x)=\Big[ \sum_{i=1}^n i \cos \Big( i + (i+1) x \Big) \Big] \cdot
\Big[ \sum_{i=1}^n i \cos \Big(i + (i+1) y \Big) \Big], \ee
which has 18 global minima $f_* \= -186.7309$ for $n=5$ in the search domain $-10 \le x,y \le 10$. \\

\begin{table}[ht]

\caption{Comparison of algorithm performance  in terms of number of iterations. \label{table-fpa}}
\centering

\begin{tabular}{rrrrr}
\hline \hline
Functions/Algorithms & GA & PSO & FPA \\
\hline

 Michalewicz ($d\!\!=\!\!16$)  & $89325 \pm 7914 (95 \%)$  & $6922 \pm 537 (98\%)$  & $3341 \pm 649 (100\%)$ \\

 Rosenbrock ($d\!\!=\!\!16$) & $55723 \pm 8901 (90\%)$ & $32756 \pm 5325 (98\%)$ & $5532 \pm 1464 (100\%) $ \\

 De Jong ($d\!\!=\!\!256$) & $25412 \pm 1237 (100\%)$ & $17040 \pm 1123 (100\%)$ & $4245 \pm 545 (100\%)$\\
 Schwefel ($d\!\!=\!\!128$) & $227329 \pm 7572 (95\%)$ & $14522 \pm 1275 (97\%)$ & $6851 \pm 448 (100\%)$ \\

 Ackley ($d\!\!=\!\!128$) & $32720 \pm 3327 (90\%)$ & $23407 \pm 4325 (92\%)$ & $3357 \pm 968 (100\%)$ \\

 Rastrigin & $110523 \pm 5199 (77 \%)$ & $79491 \pm 3715 (90\%)$ & $10840 \pm 2689 (100\%)$ \\

 Easom & $19239 \pm 3307 (92\%)$ & $17273 \pm 2929 (90\%)$ & $4017 \pm 982 (100\%)$ \\

 Griewank & $70925 \pm 7652 (90\%)$ & $55970 \pm 4223 (92\%)$ & $4918 \pm 1429 (100\%)$ \\

 Yang ($d=16$) & $27923 \pm 3025 (83\%)$ & $14116 \pm 2949 (90\%)$ & $4254 \pm 1839 (100\%)$ \\

Shubert(18 minima) & $54077 \pm 4997 (89\%)$ & $23992 \pm 3755 (92\%)$ & $9271 \pm 1758 (100\%)$ \\

\hline
\end{tabular}
\end{table}

For the above ten test functions, each function can have varied dimensions; then there is an issue
which dimensions should be used in the simulations. Research suggests that higher-dimensional problems
tend to be more challenging, and a new algorithm should be tested against a wide range of
functions in terms of function properties and dimensions. Therefore, we tend to focus on
problems with higher dimensions.

In addition, we have used three algorithms to find their optimal solution
with a given tolerance $10^{-5}$. The three algorithms are genetic algorithm (GA), particle swarm
optimization (PSO) and the new flower pollination algorithm (FPA). For each algorithm, we have
carried out $100$ independent runs using a population size $n=25$ and $p=0.8$ for FPA,
crossover probability $0.95$ and mutation probability $0.05$ for GA, and learning parameters $2$ for PSO.
The results are summarized in Table \ref{table-fpa}.
In the table, the results are provided as mean $\pm$ standard deviation (success rate).
For example, $3341 \pm 649 (100\%)$ means that mean number iterations is 3341 with one standard
deviation of $649$ and a success rate of $100\%$. The total number of function evaluations is $n$ times the mean number of iterations. For example,
the number of iterations is 3341 in the table, so the total number of function evaluations
is $3341n=3341 \times 25=83525$.

\subsection{Design Optimization}

Pressure vessels are literally everywhere such as champagne bottles and gas tanks.
For a given volume and working pressure, the basic aim of designing
a cylindrical vessel is to minimize the total cost. Typically, the design variables are
the thickness $d_1$ of the head, the thickness $d_2$ of the body, the inner radius
$r$, and the length $L$ of the cylindrical section \cite{Cag}.
This is a well-known test problem for optimization  and it can be written as
\be \textrm{minimize } f(\x) = 0.6224 d_1 r L + 1.7781 d_2 r^2
 + 3.1661 d_1^2 L + 19.84 d_1^2 r, \ee
subject to the following constraints
\be
\begin{array}{lll}
 g_1(\x) = -d_1 + 0.0193 r \le 0 \\
 g_2(\x) = -d_2 + 0.00954 r \le 0 \\
 g_3(\x) = - \pi r^2 L -\frac{4 \pi}{3} r^3 + 1296000 \le 0 \\
 g_4(\x) =L -240 \le 0.
\end{array}
\ee

The simple bounds are  \be 0.0625 \le d_1, d_2 \le 99 \times 0.0625, \ee
and \be 10.0 \le r, \quad L \le 200.0. \ee

Recently, Cagnina et al (2008) used an efficient particle swarm optimiser
to solve this problem and they found the best solution
\be f_* \approx 6059.714, \ee
at \be \x_* \approx (0.8125, \; 0.4375, \; 42.0984, \; 176.6366). \ee
This means the lowest price is about $\$6059.71$.

Using the proposed flower pollination algorithm, we have easily found the same as the solution
$f_* \approx 6059.714$ obtained by Cagnina et al \cite{Cag,Gandomi}.

\begin{figure}
\centerline{\includegraphics[width=3in,height=2.5in]{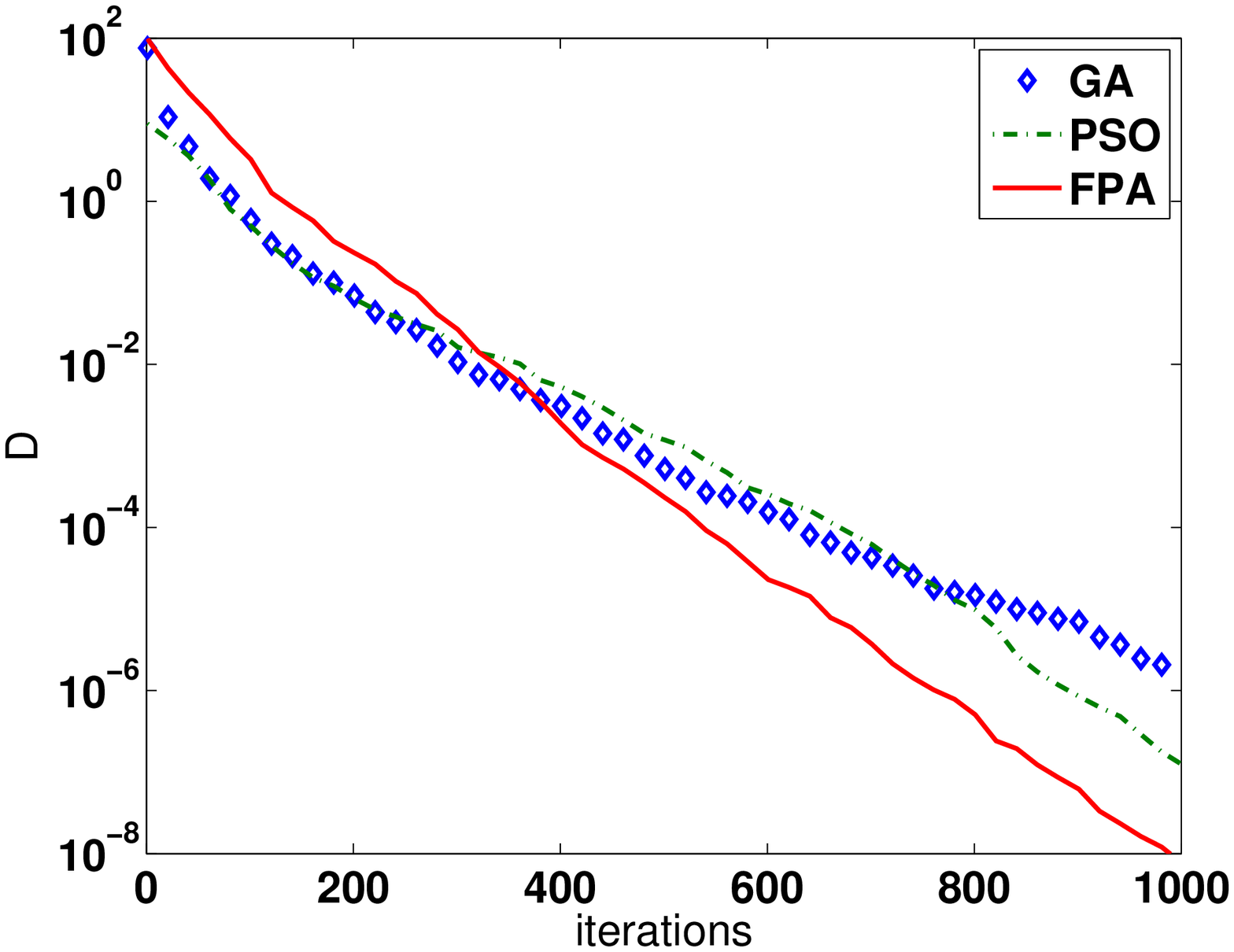} }
\caption{Error variations and comparison of GA, PSO and FPA. \label{pafig-40} }
\end{figure}

The current best solution can be stored during iterations. Then, we can calculate the
errors $D$ in terms of the difference between the current solution to the best mean solution
after $40$ independent runs. Figure \ref{pafig-40} shows the typical variations of $D$ during iterations.
We can also see that the proposed algorithm approaches the optimal solution exponentially (see  Fig. \ref{pafig-40}).

Among the three methods, the proposed FPA obtained the best result and converged most quickly.

\section{Discussions}

Flowering plants have evolved some interesting features of flower pollination, and we have successfully
developed a new flower algorithm to mimic these characteristics. Our simulation results have shown that
the the proposed flower pollination algorithm is very efficient and can outperform both genetic algorithm
and particle swarm optimization. The convergence rate is essentially exponential as we have seen from the
convergence comparison in the previous section.

The reasons that FPA is efficient can be twofold: long-distance pollinators and flower consistency.
Pollinators such as insects can travel long distance, and thus they introduce the ability (into the
algorithm) that they can escape any local landscape and subsequently explore larger search space.
This acts as exploration moves. On the other hand, flower consistency ensure
that the same species of the flowers (thus similar solutions) are chosen
more frequently and thus guarantee the convergence more quickly. This step is essentially
an exploitation step. The interplay and interaction of these key components and the selection of the
best solution $\ff{g}_*$ ensure that the algorithm is very efficient.

\section{Conclusions}

In the present algorithm, for simplicity, we have assumed that each flower only produce one pollen
gamete, this simplifies the implementation greatly. However, to assign each flower with multiple pollen gametes
and each plant with multiple flowers can have some advantages for some applications such as
image compression, multiobjective optimization, and graph colouring. This can form a topic
for further research.

Another possible extension is to use design other schemes for flower constancy. At present,
this constancy is realized by a simple formula. Other more exotic form may be useful to
certain type of problem, though the exact improvement may need some extensive simulations.

Furthermore, it is possible to extend the flower algorithm to a discrete version so that it can
solve combinatorial optimization problems. All these extensions will be very useful. We hope that
this paper will inspire more active research in metaheuristics in the near future.

%% End of text %%

\end{document}